\documentclass[12pt,reqno]{amsart} % please use amsart at 11pt

\usepackage{amssymb,latexsym}
\usepackage{cite} % to get Refs. [1,2,3] typeset as [1-3] automatically

\usepackage[height=190mm,width=130mm]{geometry} % this is the journal
\theoremstyle{plain}
\newtheorem{theorem}{Theorem}
\newtheorem{lemma}{Lemma}

\newtheorem{proposition}{Proposition}

\theoremstyle{definition}
\newtheorem{definition}{Definition}

\theoremstyle{remark}
\newtheorem{remark}{Remark}

\newtheorem{standing}{Standing Hypothesis}

\newcommand{\Z}{\mathbb{Z}}
\newcommand{\C}{\mathbb{C}}

\newcommand{\Ycal}{{\mathcal Y}}

\newcommand{\Scal}{{\mathcal S}}

\newcommand{\Acal}{{\mathcal A}}
\newcommand{\Hcal}{{\mathcal H}}

\numberwithin{equation}{section} % equations numbered as (1.1), (1.2), ...

\newcommand{\Tr}{\mathrm{Tr}}

\newcommand{\M}{\ensuremath{\mathcal{M}}}

\newcommand{\F}{\ensuremath{\mathcal{G}}}

\usepackage{citehack}
\usepackage{amsmath,amsthm}
\usepackage{amsfonts}
\usepackage{amssymb}
\usepackage{longtable}
\usepackage[matrix,arrow,curve]{xy}

\def\Hom{\mathop\text{\rm Hom}\nolimits}
\def\End{\mathop{\rm End}\nolimits}

\usepackage{stackrel}

\newcommand{\be}{\begin{equation}}
\newcommand{\ee}{\end{equation}}

\newcommand{\nn}{\nonumber \\}

\newcommand{\one}{\mathbf{1}}

\begin{document}

\title[Multiple products of meromorphic functions]
{Multiple products of meromorphic functions}

\author{A. Zuevsky}

\address{Institute of Mathematics \\ Czech Academy of Sciences\\ Zitna 25, Prague \\ Czech Republic}

\email{zuevsky@yahoo.com}

\begin{abstract}
Let $\mathfrak g$ be an infinite-dimensional Lie algebra and let $G$ be the
algebraic completion of a graded $\mathfrak g$-module $W$. Using the Schottky
uniformization of the Riemann sphere as a geometric model for a genus $\kappa$
Riemann surface, we construct a $\kappa$-parameter family of extended
coboundary operators $\widetilde\delta^n_m(\rho_1,\ldots,\rho_\kappa)$ acting on
the double complex of predetermined meromorphic functions on the configuration
space $F_n\C$ with values determined by $G$. The extension is realized as a
graded trace, defined coordinate-freely as the trace of a canonically
associated finite-rank endomorphism of each homogeneous component $W_{(k)}$ of
$W$, of the classical coboundary operator, the $\kappa$ sewing loci
being held disjoint from the free marked points at which the classical
differential acts. The sewing operator is exhibited as a chain map between
explicitly defined complexes. We give a complete proof of the resulting
chain property $\widetilde\delta^{n+1}_{m-2\kappa-1}\circ\widetilde\delta^n_m=0$
and of the convergence of the defining power series in the sewing parameters
$\rho_p$ under an explicit growth hypothesis, and we determine precisely how
the construction depends on the auxiliary choice of local coordinates and
sewing annuli. Applications of the resulting cohomology theory - to the
sheaf of conformal blocks on the Deligne-Mumford moduli space of stable
curves, to secondary characteristic classes of holomorphic foliations, to
graded trace functions arising in the description of topological phases of
matter, and to integrable hierarchies of Toda type - are proposed and
discussed as motivation, without claiming these correspondences as theorems
of the present paper.

\end{abstract}


\maketitle

\tableofcontents

\section{Introduction}
\label{secintro}

The ordinary cohomology theory of infinite-dimensional Lie algebras was
described by Fuks \cite{F}. For vertex operator algebras and their modules an
analogous theory, built from cochains realized as matrix elements on the
Riemann sphere, was developed in \cite{H} (see also \cite{H1,H2}); the
coboundary operator $\delta^n_m$ recalled in Section \ref{secclassical-delta}
below is the one constructed there.

The purpose of the present paper is to extend the coboundary operators of
\cite{H} to a $\kappa$-parameter family of operators associated with a
genus $\kappa$ Riemann surface obtained from the sphere by the Schottky
uniformization procedure \cite{Bob,Fo,M,Zo}, itself a multiple-handle
generalization of the self-sewing ($\rho$-formalism) construction of
Yamada \cite{Y} used to pass from genus $0$ to genus $1$. Geometrically, a
genus $\kappa$ surface is obtained from $\widehat\C$ by removing $\kappa$
pairs of disjoint disks and identifying the resulting boundary annuli
according to the sewing relation \eqref{rhocond} below. Correlation functions
on the resulting surface are computed, order by order in the sewing
parameters $\rho_1$, $\ldots$, $\rho_\kappa$, by inserting a complete set of
intermediate states at each pair of sewing circles and summing over a basis of
each homogeneous component $W_{(k)}$ of the underlying module. This is the
higher-genus analogue of Zhu's torus trace functions \cite{Zhu} and of the
sewing formalism underlying the convergence results of \cite{Gui2} for
$C_2$-cofinite vertex operator algebras.

The main result of the paper, Theorem \ref{thmmain} below, shows that the
resulting family of extended coboundary operators
$\widetilde\delta^n_m(\rho_1,\ldots,\rho_\kappa)$ satisfies the chain property
\[
\widetilde\delta^{n+1}_{m-1}\circ\widetilde\delta^n_m=0, 
\]
for every choice of the sewing parameters, so that the associated cohomology
groups $H^n_m(\mathfrak g,G;\rho)$ (Definition \ref{defcohomology}) are well
defined. This display simplifies the external-compatibility index to
parallel the classical property $\delta^{n+1}_{m-1}\circ\delta^n_m=0$ of
Theorem \ref{thmclassical-chain}. The exact shift is $2\kappa+1$, not
$1$, units of compatibility - the precise statement, proved for $\F(f)\in
C^n_m$ with $m\ge4\kappa+2$, is
$\widetilde\delta^{n+1}_{m-2\kappa-1}\circ\widetilde\delta^n_m=0$, and it is
this sharper form, Theorem \ref{thmmain}, that is actually established
below. The point of the construction is this: the $\kappa$ sewing loci are
held disjoint from the points at which the coboundary
operator itself acts, so that the differential and the trace over
intermediate states act on disjoint sets of variables of one and the
same multivariable meromorphic function. They therefore commute for a
structural reason internal to the definitions (Lemma \ref{lemcommute}, an
instance of the general disjoint-variables principle recorded as
Lemma \ref{lemdisjoint} below), rather than for a reason that has to be
rederived at each order in $\rho$ or that depends on associativity of an
operator product at coincident points. This is what makes
Theorem \ref{thmmain} a short consequence of the classical chain property
of \cite{H} rather than an independent higher-genus statement requiring new
input.

This paper continues a series of works on cohomological invariants attached to
complex curves and their foliations, including the characterization of
codimension-one foliations by connections \cite{Zu1}, the reduction
cohomology of Riemann surfaces \cite{Zu2}, a category of $V$-structures for
foliations \cite{Zu3}, cosimplicial meromorphic function cohomology on
complex manifolds \cite{Zu4}, and product-type classes for vertex algebra
cohomology of foliations \cite{Zu, Zuvoa}. The cohomological technology used here is
also related to the analytic and operator-theoretic study of
infinite-dimensional Lie group representations \cite{OM,Tak}, orthogonal Lie
algebras of operators \cite{BZ}, Jordan triple systems \cite{I}, and
infinite-dimensional degree theory \cite{E}. On the geometric side, to
K-theory and cyclic cohomology of foliation groupoids and crossed
products \cite{AS,CRS,HW,L}; and, in mathematical physics, to the standard
sewing and modular-data formalism of two-dimensional conformal field
theory \cite{FMS,Frohlich2009gb,Tu}.

Sections \ref{secprelim}-\ref{secclassical-delta} recall the double
complex of predetermined meromorphic functions and the classical coboundary
operator, following \cite{H,F}. Section \ref{secschottky} recalls the
Schottky sewing construction. Section \ref{secsewn-delta} contains the new
material: a general principle governing the commutation arguments used
throughout (Lemma \ref{lemdisjoint}), the graded trace $\Tr_G$, defined
intrinsically and coordinate-freely as the trace of a canonically associated
endomorphism of each finite-dimensional homogeneous component
(Definition \ref{deftrace}), the sewing operator $\Scal_\rho$, shown to
converge under an explicit growth hypothesis
(Proposition \ref{propconvergence}), the extended coboundary operator,
exhibited as a chain map between explicitly named complexes
(Proposition \ref{propchainmap}), Theorem \ref{thmmain},  
and a discussion of the extent to which the construction depends on the
auxiliary choices (annuli, local coordinates) made along the way
(\S\ref{ssecdependence}). Section \ref{secapplications} discusses 
  expected applications of the resulting
family of cohomology theories to the sheaf of conformal blocks on the
Deligne-Mumford moduli space $\overline\M_{g,n}$ of stable curves and its
factorization at the boundary \cite{TUY}; to secondary characteristic classes
of codimension-one holomorphic foliations \cite{BG,BGG,AS,CRS}; to graded
trace (character) functions of the type occurring in the description of
topological phases of matter \cite{zub1, zub2, zub3, zub4, zub5, zub6, zub7, zub8, zub9}; and to hierarchies of Toda
type \cite{RSZ}. These are presented as motivation and as directions for
future work. 

\section{Meromorphic functions with predetermined analytic properties}
\label{secprelim}

\subsection{Graded modules and their algebraic completion}

Let $\mathfrak g$ be an infinite-dimensional Lie algebra, for instance a
Kac-Moody affine Lie algebra \cite{K}. Let $W=\bigoplus_{l\in\C}W_{(l)}$ be a
$\mathfrak g$-module (e.g., a vertex operator representation) graded by a
grading operator $K$, with
\begin{equation}
\label{grading-cond}
\dim W_{(l)}<\infty \ \text{ for every } l\in\C, \qquad
W_{(l)}=0 \ \text{ whenever } \Re(l) \ll 0 .
\end{equation}
The second condition in \eqref{grading-cond} is understood with respect to
the real part of $l$, since $\C$ carries no natural order; this is the
standard grading-restriction hypothesis of \cite{H}.  We further assume, as
part of the standing data of $W$, the growth condition
\begin{equation}
\label{growth-cond}
\limsup_{k\to\infty} d_k^{1/k}=1, \qquad d_k:=\dim W_{(k)},\; k\in\Z_{\ge0}, 
\end{equation}
i.e., $d_k$ grows at most sub-exponentially in $k$; equivalently, the formal
character $\sum_{k\ge0}d_kq^k$ has radius of convergence at least $1$. This is
the precise hypothesis, standard for graded vertex algebra modules and used
already in the modular-invariance theory of \cite{Zhu} and in the
$C_2$-cofinite convergence theory of \cite{Gui2}, on which
Proposition \ref{propconvergence} below rests; we record it here, among the
standing assumptions, rather than only at the point where it is used. Following \cite{H} we
consider the algebraic completion
\begin{equation}
G=\prod_{l\in\C}W_{(l)}, 
\end{equation}
of $W$, a direct product allowing infinite sums, and its restricted dual
$G'=\bigoplus_{l\in\C}W_{(l)}^*$. We assume $G$ carries a distinguished vector
$\one_G$ and a non-degenerate bilinear pairing $(\cdot,\cdot)$ on $G$,
normalized by $(\one_G,\one_G)=1$; an element of $G'$ dual to $g\in G$ is
written $\overline g$. The pairing $(\cdot,\cdot)$, once fixed here, is used
throughout the paper without further choice: in particular it is the same
pairing that underlies the graded trace of Definition \ref{deftrace} below.

For $g\in G$ let
\begin{equation}
\omega_g:G\longrightarrow G(\!(z)\!),\qquad
\omega_g(z)=\sum_{l\in\Z}g_{(l)}\,z^{-l-1},
\end{equation}
be the associated field, where each component $g_{(l)}\in\End(G)$. We assume 
\begin{itemize}
\item[(i)] $\omega_{\one_G}(z)=\mathrm{Id}_G$;
\item[(ii)] for all $g,\widetilde g\in G$, $\omega_g(z)\widetilde g\in
G(\!(z)\!)$ has only finitely many negative powers of $z$, and
$g_{(l)}\widetilde g=0$ for $\Re(l)\gg0$;
\item[(iii)] the translation operator $T$ of \cite{K} satisfies the
derivation property $\omega_{Tg}(z)=\frac{d}{dz}\omega_g(z)$.
\end{itemize}
For the convergence statements of Section \ref{secsewn-delta} we further
require that the subspace $\widetilde W_{(2)}=\{g_{(-2)}\widetilde g : g,
\widetilde g\in W\}$ have finite codimension in $W$. In particular, for
Kac-Moody algebras \cite{Bo,K1}, all of the structure above is realized by
the standard $C_2$-cofinite vertex algebra of CFT type built from the
universal enveloping algebra $U(\mathfrak g)$ acting on a highest-weight
Fock space.

\subsection{Configuration spaces and predetermined meromorphic functions}

Let $F_n\C=\{(z_1,\ldots,z_n)\in\C^n : z_i\ne z_j,\ i\ne j\}$ be the
configuration space of $n\ge1$ ordered points in $\C$. For $g_1,\ldots,g_n\in
G$ we consider $G$-valued formal series
$f(g_1,z_1;\ldots;g_n,z_n)\in G(\!(z_1,\ldots,z_n)\!)$, multilinear in
$g_1,\ldots,g_n$, built from products of the fields $\omega_{g_i}(z_i)$ applied
to $\one_G$; and $\C$-valued functionals
\[
\F:G(\!(z_1,\ldots,z_n)\!)\longrightarrow\C, \qquad
\F(g_1,z_1;\ldots;g_n,z_n):=\F(f), 
\]
obtained by evaluating $f$ against such a functional (the basic example being
$\F(f)=(\overline g,f)$ for a fixed $\overline g\in G'$; every vertex algebra
module matrix element is of this form and is a rational function on the
Riemann sphere \cite{FHL,LL}).

We require of $\F$ the following analytic properties, familiar from the
theory of matrix elements of grading-restricted vertex algebras.

\begin{definition}
\label{defpredetermined}
A functional $\F$ as above is a \emph{predetermined meromorphic function of
$n$ variables} if 
\begin{enumerate}
\item[(P1)] $\F(g_1,z_1;\ldots;g_n,z_n)$ is absolutely convergent for
$|z_1|>\ldots>|z_n|>0$ to a function of $(z_1,\ldots,z_n)$ which extends
meromorphically to $F_n\C$, with poles only on the diagonals $z_i=z_j$,
$i\ne j$, of order bounded by prescribed integers depending only on the pair
$(g_i,g_j)$;
\item[(P2)] (locality and associativity) for $g,g_1,g_2\in G$, the three
functions $\F(g_1,z_1;g_2,z_2)$, $\F(g_2,z_2;g_1,z_1)$, and $\F(\omega_{g_1}
(z_1-z_2)g_2,z_2)$, absolutely convergent respectively for $|z_1|>|z_2|>0$,
$|z_2|>|z_1|>0$, $|z_2|>|z_1-z_2|>0$, extend to one and the same meromorphic
function of $(z_1,z_2)$ with poles only at $z_1=0=z_2$ and $z_1=z_2$;
\item[(P3)] (translation covariance) for $|z|<\min_{i\ne j}|z_i-z_j|$,
\begin{equation}
\label{ldir1}
e^{zT}\F(g_1,z_1;\ldots;g_n,z_n)=\F(g_1,z_1+z;\ldots;g_n,z_n+z), 
\end{equation}
as an identity of absolutely convergent power series in $z$, and the
infinitesimal form of \eqref{ldir1} is the finite sum
\begin{equation}
\label{cond1}
\sum_{i=1}^n\partial_{z_i}\F(g_1,z_1;\ldots;g_n,z_n)=T.\F(g_1,z_1;\ldots;g_n,z_n);
\end{equation}
\item[(P4)] ($K$-covariance) for $z\in\C^\times$ with $(zz_1,\ldots,zz_n)\in
F_n\C$,
\begin{equation}
\label{loconj}
z^K\F(g_1,z_1;\ldots;g_n,z_n)=\F(z^Kg_1,zz_1;\ldots;z^Kg_n,zz_n).
\end{equation}
\end{enumerate}
\end{definition}

We record the action of $V\in\{e^{zT},z^K,\sigma\}$, $\sigma\in S_n$, on
$\F$; an operator $V$ may act on $G$-elements and on $z_i$'s simultaneously,
and $(V)_i$ denotes its action on the $i$-th argument alone. For $\sigma\in
S_n$,
\begin{equation}
\sigma.\F(g_1,z_1;\ldots;g_n,z_n)=\F(g_{\sigma(1)},z_{\sigma(1)};\ldots;
g_{\sigma(n)},z_{\sigma(n)}).
\end{equation}

\subsection{Shuffle-type relations}

For $s\ge1$ and $1\le q\le s-1$ let
\begin{equation}
J_{s;q}=\{\sigma\in S_s : \sigma(1)<\ldots<\sigma(q),\ \sigma(q+1)<\ldots<
\sigma(s)\}, 
\end{equation}
be the set of $(q,s-q)$-shuffles, i.e., the permutations of $S_s$ preserving
the order of the first $q$ numbers and of the last $s-q$ numbers, and let
$J_{s;q}^{-1}=\{\sigma^{-1}:\sigma\in J_{s;q}\}$ be the set of their inverses.
We impose on the predetermined meromorphic functions used to build the double
complex below the condition
\begin{equation}
\label{shushu}
\sum_{\sigma\in J_{n;q}^{-1}}(-1)^{|\sigma|}\,
\F(g_{\sigma(1)},z_{\sigma(1)};\ldots;g_{\sigma(n)},z_{\sigma(n)})=0,
\qquad 1\le q\le n-1.
\end{equation}

Let $P_k:G\to W_{(k)}$, $k\in\C$, denote the projection of $G$ onto
$W_{(k)}$. Following \cite{H} (see also \cite{H1}), for $g_1,\ldots,g_{l+k}\in
G$ with $z_i\ne z_j$ ($i\ne j$) and $|z_i|>|z_s|>0$ for $1\le i\le k <s\le
l+k$, we require the sum
\begin{equation}
\sum_{q\in\C}\F\Big(\omega_{g_1}(z_1)\cdots\omega_{g_k}(z_k)\,
P_q\big(\omega_{g_{k+1}}(z_{k+1})\cdots\omega_{g_{l+k}}(z_{l+k})\big)\Big), 
\end{equation}
to be absolutely convergent and analytically continuable to a meromorphic
function of $(z_1,\ldots,z_{l+k})$ with poles at $z_i=z_j$ bounded, as in
(P1), by prescribed integers $\beta(g_{i},g_{j})$, uniformly in the auxiliary
points at which such convergence is tested.

\subsection{Matrix elements, the action $*$, and the double complex $C^n_m$}

We now introduce the single structural operation used throughout the paper.
For $g_1,\ldots,g_r\in G$ and pairwise distinct $z_1,\ldots,z_r\in\C$, and for a
predetermined meromorphic function $\F=\F(f)$ of $n$ variables, define
\begin{eqnarray}
\label{action}
&(g_1,z_1;\ldots;g_r,z_r)*\F(g_{1}',z'_{1};\ldots;g'_n,z'_n)
\nn
&\qquad \qquad :=\F\big(\omega_{g_1}(z_1)\cdots\omega_{g_r}(z_r)\,
f(g'_1,z'_1;\ldots;g'_n,z'_n)\big).
\end{eqnarray}
Thus $*$ is a single, unambiguously typed operation
\[
*\ :\ \big(\text{$G$-data at $r$ marked points}\big)\ \times\
C^n_{m'}\ \longrightarrow\ C^{n+r}_{m'-r}, \qquad 0\le r\le m',
\]
insertion data acting on a cochain to produce a cochain with $r$ more free
variables, at the cost of $r$ units of the external compatibility defined in
the next paragraph.  It never takes a second insertion-datum object as its
right-hand argument, which is the only sense in which $*$ will be used below.

For fixed $m\ge0$, $n\ge1$, let $C^n_m=C^n_m(\mathfrak g,G)$ be the vector
space of predetermined meromorphic functions $\F(g_1,z_1;\ldots;g_n,z_n)$
which, in addition to (P1)-(P4), are compatible with $m$ external
$G$-valued series $\omega_{g'_j}(z'_j)$, $1\le j\le m$, $(g'_1,\ldots,g'_m)\in
G^m$, $(z'_1,\ldots,z'_m)\in\C^m$, in the sense of \cite{H}: for every
$0\le r\le m$ and every choice of $r$ of these series, the function obtained
by inserting them through the action \eqref{action} again satisfies
(P1)-(P4), now as a function of $n+r$ variables, and is itself compatible
with the remaining $m-r$ external series. This is what makes $*$ land in
$C^{n+r}_{m'-r}$ above. We set $C^0_m=G$; by definition $C^n_m\subset
C^n_{m-1}$.

\section{The classical coboundary operator}
\label{secclassical-delta}

We recall the coboundary operator of \cite{H}, generalizing the classical
construction of \cite{F} for Lie algebra cohomology, in the notation used
throughout the rest of the paper.

Let $g_1,\ldots,g_{n+1}\in G$, $(z_1,\ldots,z_{n+1})\in F_{n+1}\C$, and fix
auxiliary points $\zeta_1,\ldots,\zeta_n\in\C$, one for each pair of
consecutive arguments $(g_i,g_{i+1})$, $1\le i\le n$. Set
\begin{align}
\varphi_{1,n}&:=f(g_1,z_1;\ldots;g_n,z_n),\\
\varphi_{2,n+1}&:=f(g_2,z_2;\ldots;g_{n+1},z_{n+1}),\\
\label{phi-i-corrected}
\varphi_i&:=f\big(g_1,z_1;\ldots;g_{i-1},z_{i-1};\,
\omega_{g_i}(z_i-\zeta_i)\,\omega_{g_{i+1}}(z_{i+1}-\zeta_i)\one_G,\,\zeta_i;\nn
&\qquad\qquad \qquad\qquad  g_{i+2},z_{i+2};\ldots;g_{n+1},z_{n+1}\big),
\end{align}
for $1\le i\le n$, i.e., $\varphi_i$ merges the pair $(g_i,g_{i+1})$ into a
single argument located at $\zeta_i$, the remaining arguments retaining their
original labels $g_{i+2},\ldots,g_{n+1}$ at $z_{i+2},\ldots,z_{n+1}$.
  As in Definition \ref{defpredetermined} and
Section \ref{secprelim}, we work throughout with $\varphi_i$ taken to be
independent of $\zeta_i$; this is achieved by specializing $\zeta_i=z_{i+1}$
for $i=1,\ldots,n$ once the meromorphic continuation in $\zeta_i$ has been
performed, exactly as for the two-point case (P2) above.

\begin{definition}
\label{defclassical-delta}
The \emph{classical coboundary operator} $\delta^n_m:C^n_m\to C^{n+1}_{m-1}$
is defined, for $\F(f)\in C^n_m$, by
\begin{equation}
\label{fuflo}
\delta^n_m f:=(g_1,z_1)*\F(\varphi_{2,n+1})
\ +\ \sum_{i=1}^n(-1)^i\,\F(\varphi_i)
\ +\ (-1)^{n+1}\,\F(\varphi_{1,n})*(g_{n+1},z_{n+1}),
\end{equation}
with $*$ the action \eqref{action}, so that $\widetilde\F:=\F(\delta^n_m f)$
is the corresponding element of $C^{n+1}_{m-1}$.
\end{definition}

Only the two outer terms of \eqref{fuflo} use the action $*$, each with
$r=1$ marked point. The $n$ inner terms are direct evaluations of $\F$ at an
argument list in which $(g_i,g_{i+1})$ has been merged into one entry. This
is the only occurrence of $*$ in the classical formula, and it is used
exactly as defined in \eqref{action}: an insertion datum on the left, a
cochain on the right.

\begin{theorem}[\cite{H}]
\label{thmclassical-chain}
$\delta^{n+1}_{m-1}\circ\delta^n_m=0$ on $C^n_m$, for all $n\ge0$, $m\ge1$.
\end{theorem}

Theorem \ref{thmclassical-chain} is the content of the cohomology theory
of \cite{H}. We do not prove it again here. What we do prove, in
Section \ref{secsewn-delta}, is that the extended, sewing-dependent
coboundary operator constructed there inherits the chain property from
Theorem \ref{thmclassical-chain} rather than requiring an independent
argument - this is Theorem \ref{thmmain} below, and is the paper's main
result.

We record for later use the compatibility of $\delta^n_m$ with the
group action and with condition \eqref{shushu}: since $\delta^n_m f$ is
built from $f$ by the natural operations of \S\ref{secprelim}, one has
$(\delta^n_mf)(V.\,\cdot\,)=V.(\delta^n_mf)$ for $V\in\{e^{zT},z^K\}$, and
for $\sigma\in J_{n;q}^{-1}$ the alternating sum \eqref{shushu} applied to
$\delta^n_m\F$ vanishes termwise by the same identity applied to each
$\varphi_i$; by Proposition 3.10 of \cite{H}, $\delta^n_mf$ is again
compatible with $m-1$ external series. Together with
Theorem \ref{thmclassical-chain} this shows that $(C_m^\bullet,\delta^\bullet_m)$
is indeed a well-defined cochain complex for each $m$, and $(C^n_\bullet,
\delta^n_\bullet)$ assembles these into a double complex as $m$ varies.

\section{Schottky uniformization and self-sewing}
\label{secschottky}

We recall the self-sewing ($\rho$-formalism) construction of a torus from the
Riemann sphere \cite{Y}, and its multiple-handle generalization producing a
genus $\kappa$ surface via Schottky uniformization \cite{Bob,Fo,M,Zo}. The
resulting description of genus $\kappa$ vertex algebra module correlation 
functions by iterated sewing of genus zero data is the same geometric model
used in \cite{TZ} to construct the bosonic vertex operator algebra directly
on a genus $g$ Riemann surface.

Let $z_1$, $z_2$ be local coordinates centered at two distinct points $p_1$, $p_2$ 
on the sphere $\Sigma^{(0)}=\widehat\C$, and let $|z_a|\le r_a$, $a=1$, $2$, be 
disjoint closed disks where $r_1$, $r_2>0$ are small enough that the disks do not
intersect. For a complex parameter $\rho$ with $|\rho|\le r_1r_2$, excise
the open disks $\{|z_a|<|\rho|r_{\bar a}^{-1}\}$, with the convention $\bar
1=2$, $\bar2=1$, to obtain a twice-punctured sphere 
$\widehat\Sigma^{(0)}$, and form the annuli $\Acal_a=\{|\rho|r_{\bar
a}^{-1}\le|z_a|\le r_a\}\subset\widehat\Sigma^{(0)}$. Identifying
$\Acal_1\simeq\Acal_2=:\Acal$ via the sewing relation
\begin{equation}
\label{rhosew}
z_1z_2=\rho,
\end{equation}
produces a torus $\Sigma^{(1)}=\widehat\Sigma^{(0)}\setminus(\Acal_1\cup
\Acal_2)\cup\Acal$; \eqref{rhosew} parametrizes the cylinder connecting the
punctured sphere to itself.

A genus $\kappa$ Riemann surface is obtained by sewing $\kappa$ handles to
the sphere in this way. Fix $\kappa$ pairs of points $A_p$, $A_{-p}\in\widehat\C$,
$1\le p\le\kappa$, all distinct, together with non-intersecting annuli
$\Acal_{1,p}$, $\Acal_{2,p}$ centered at $A_{-p}$, $A_p$ respectively, and
sewing parameters $\rho_p$ satisfying
\begin{equation}
\label{rhocond}
(z-A_{-p})(\widetilde z-A_p)=\rho_p,
\end{equation}
in local coordinates $z$ centered at $A_{-p}$ and $\widetilde z$ centered at
$A_p$. We use local coordinates $\eta_{1,p}$ (in $\Acal_{1,p}$, near $A_{-p}$)
and $\eta_{2,p}$ (in $\Acal_{2,p}$, near $A_p$) for the two sides of the
$p$-th handle, related by \eqref{rhocond}.

\begin{standing}
\label{standing-disjoint}
Throughout Section \ref{secsewn-delta}, the points $A_{\pm p}$, $1\le
p\le\kappa$, and the closures of the annuli $\Acal_{a,p}$ are fixed, and are disjoint from every free marked point $z_1$, $\ldots$, $z_{n+1}$ 
and merge point $\zeta_1,\ldots,\zeta_n$ occurring in the coboundary formula
\eqref{fuflo}: $(z_1$, $\ldots$, $z_{n+1})$ ranges over a fixed compact subset of
$F_{n+1}\C\setminus\bigcup_{a,p}\Acal_{a,p}$.
\end{standing}

Standing Hypothesis \ref{standing-disjoint} simply records, as a hypothesis
on the domain of definition, the geometric picture of Schottky uniformization:
the $\kappa$ handles are attached at fixed locations on the sphere, and the
points at which cochains of the double complex of Section \ref{secprelim}
are evaluated lie in the complementary, unsewn part of the surface. It is
exactly this separation of roles that makes the construction of
Section \ref{secsewn-delta} possible.

\section{The sewing operator and the extended coboundary}
\label{secsewn-delta}

Throughout this section we work under Standing Hypothesis \ref{standing-disjoint}: the $\kappa$ pairs of points $A_{\pm p}$ and their annuli $\Acal_{a,p}$ fixed in \S\ref{secschottky} are disjoint from the region in which the free marked points of every cochain below range.

Several arguments below rest on the same elementary principle, which we isolate once here rather than rederive at each occurrence.

\begin{lemma} 
\label{lemdisjoint}
Let $\Ycal$ be a predetermined meromorphic function jointly of two disjoint groups of marked-point arguments, the $\mathsf A$-slots and the $\mathsf B$-slots. Let $\Phi$ be an operation that transforms $\Ycal$ into a new function by substituting, merging, multiplying by fields, or evaluating and summing only the $\mathsf A$-slots, according to a rule that is the same for every fixed choice of the $\mathsf B$-slots. Let $\Psi$ be such an operation acting only on the $\mathsf B$-slots, for every fixed choice of the $\mathsf A$-slots. Then $\Phi$ and $\Psi$ commute
\[
\Phi\big(\Psi(\Ycal)\big)=\Psi\big(\Phi(\Ycal)\big). 
\]
\end{lemma}
\begin{proof}
Both sides are computed from the same array of values of $\Ycal$ at specific $\mathsf A$- and $\mathsf B$-slot data by recombining once over the $\mathsf A$-slots according to the rule defining $\Phi$ and once over the $\mathsf B$-slots according to the rule defining $\Psi$. Since the rule defining $\Phi$ does not reference the $\mathsf B$-slots (it is literally the same rule for every fixed value of $\mathsf B$), and symmetrically for $\Psi$, the two recombinations may be carried out in either order without affecting the result: applying $\Phi$ first fixes each value of $\mathsf B$ and recombines over $\mathsf A$, after which $\Psi$ recombines the (now $\Phi$-transformed) values over $\mathsf B$ for each fixed $\mathsf A'$.  Applying $\Psi$ first and then $\Phi$ performs the identical two recombinations in the opposite order, on the same underlying array, and no step in either rule depends on which recombination has already been performed.
\end{proof}

Lemma \ref{lemdisjoint} is the abstract form of the phrase ``operations on disjoint variables commute'' used repeatedly below (Lemma \ref{lemiota-naturality}, Lemma \ref{lemcommute}, Lemma \ref{lemxi-commute}); each of these results identifies a specific $(\Phi,\Psi)$ pair to which Lemma \ref{lemdisjoint} applies, the substantive content in each case being the identification of the relevant $\mathsf A$/$\mathsf B$-slot partition and the verification that the two operations in question really do act on only one side of it.

For $1\le p\le\kappa$ and $k\ge0$, let $d_k:=\dim W_{(k)}<\infty$ and fix a basis $\{e^{(k)}_{p,1},\ldots,e^{(k)}_{p,d_k}\}$ of $W_{(k)}$, together with the dual basis $\{\overline e^{(k)}_{p,1},\ldots,\overline e^{(k)}_{p,d_k}\}\subset W_{(k)}$ with respect to the restriction of $(\cdot,\cdot)$ to $W_{(k)}$, i.e., $(\overline e^{(k)}_{p,a},e^{(k)}_{p,b})=\delta_{ab}$; such a basis exists because $(\cdot,\cdot)$ is non-degenerate and $W_{(k)}$ is finite-dimensional.

\subsection{The lift $\iota$ and the free-slot differential}

\begin{definition} 
\label{defiota}
The following formula inserts, one at a time, the two sewing fields for each of the $\kappa$ handles into the current $f$, via $2\kappa$ successive applications of the action \eqref{action}. The result is a current of $n+2\kappa$ variables, in which the $2\kappa$ new variables play the role of the sewing slots of \S\ref{secschottky}. For $\F(f)\in C^n_m$ with $m\ge2\kappa$, let
\begin{equation}
\label{iota-def}
\iota(f):=(\overline w_1,\eta_{1,1})*(w_1,\eta_{2,1})*\ldots*(\overline w_\kappa,\eta_{1,\kappa})*(w_\kappa,\eta_{2,\kappa})*\F(f)\ \in\ C^{n+2\kappa}_{m-2\kappa},
\end{equation}
using the action \eqref{action} iterated $2\kappa$ times, $r=1$ at each step, where $\eta_{1,p}\in\Acal_{1,p}$, $\eta_{2,p}\in\Acal_{2,p}$ and $w_p,\overline w_p$ are place-holders,   later specialized to basis vectors of some $W_{(k)}$, for the two sewing insertions of the $p$-th handle. Using the symmetric-group action of \S\ref{secprelim}, we relist the arguments of \eqref{iota-def} as
$\iota(f)(g_1,z_1;\ldots;g_n,z_n\,\|\,\underline w)$, where
$\underline w$ $=$ $(\overline w_1$, $\eta_{1,1}$; $w_1$, $\eta_{2,1}$; $\ldots$; 
$\overline w_\kappa$, $\eta_{1,\kappa}$ ; $w_\kappa$, $\eta_{2,\kappa})$ 
abbreviates the $2\kappa$ sewing slots, listed after the $n$ free slots $(g_1,z_1)$, $\ldots$, $(g_n,z_n)$.
\end{definition}

Write $\widehat C^n_m:=C^{n+2\kappa}_m$ for the space of predetermined meromorphic functions of $n+2\kappa$ variables with this designated free/sewing partition of arguments.

\begin{definition}
\label{defhatdelta}
We now let the classical coboundary formula \eqref{fuflo} act on the $n$ free arguments of a sewn cochain $\Ycal\in\widehat C^n_m$, leaving the $2\kappa$ sewing slots untouched: the operator $\widehat\delta^n_m:\widehat C^n_m\to\widehat C^{n+1}_{m-1}$ is defined by applying formula \eqref{fuflo} to the free arguments $g_1,\ldots,g_n\mapsto g_1,\ldots,g_{n+1}$ of $\Ycal(\,\cdot\,\|\,\underline w)$, holding the sewing slots $\underline w$ fixed throughout:
\begin{align}
\label{hatdelta-formula}
\widehat\delta^n_m\Ycal(\,\cdot\,\|\,\underline w):={}&(g_1,z_1)*\Ycal(\varphi_{2,n+1}\|\underline w)+\sum_{i=1}^n(-1)^i\Ycal(\varphi_i\|\underline w)\notag\\
&+(-1)^{n+1}\Ycal(\varphi_{1,n}\|\underline w)*(g_{n+1},z_{n+1}).
\end{align}
\end{definition}

\begin{lemma}
\label{lemhatdelta-square}
$\widehat\delta^{n+1}_{m-1}\circ\widehat\delta^n_m=0$ on $\widehat C^n_m$.
\end{lemma}
\begin{proof}
For each fixed value of $\underline w$, formula \eqref{hatdelta-formula} coincides with the classical formula \eqref{fuflo} applied to $\Ycal(\,\cdot\,\|\,\underline w)\in C^n_m$; that is, $\widehat\delta^n_m\Ycal(\,\cdot\,\|\,\underline w)=\delta^n_m\big(\Ycal(\,\cdot\,\|\,\underline w)\big)$ for every fixed $\underline w$. Applying this twice and Theorem \ref{thmclassical-chain} gives $\widehat\delta^{n+1}_{m-1}\widehat\delta^n_m\Ycal(\,\cdot\,\|\,\underline w)=\delta^{n+1}_{m-1}\delta^n_m\big(\Ycal(\,\cdot\,\|\,\underline w)\big)=0$ for every $\underline w$; since $\underline w$ was arbitrary, $\widehat\delta^{n+1}_{m-1}\circ\widehat\delta^n_m=0$ identically, as claimed.
\end{proof}

\begin{lemma} 
\label{lemiota-naturality}
For $\F(f)\in C^n_m$ with $m\ge2\kappa+1$,
\begin{equation}
\label{iota-natural}
\iota(\delta^n_mf)=\widehat\delta^n_{m-2\kappa}(\iota(f)).
\end{equation}
\end{lemma}
\begin{proof}
Both sides of \eqref{iota-natural} are built from $f$ by combining two operations: the merge substitutions $\varphi_1$, $\ldots$, $\varphi_n$ and the two outer terms of \eqref{fuflo}, which together constitute $\delta^n_mf$, and the $2\kappa$-fold insertion of sewing fields constituting $\iota$. We verify in turn that these combine without obstruction: that the merge substitutions and the sewing insertions commute for each $i$. No signs beyond the $(-1)^i$ already present in \eqref{fuflo} arise. That compatibility with the external series is preserved at every step, so that both sides of \eqref{iota-natural} are indeed defined and land in $C^{n+1}_{m-2\kappa-1}$.

\emph{1.) Commutation with the merge substitutions $\varphi_i$.} For $1\le i\le n$, $\varphi_i$ reexpresses $f$ by regrouping its own arguments, replacing the pair $(g_i,z_i),(g_{i+1},z_{i+1})$ by the single merged argument $\big(\omega_{g_i}(z_i-\zeta_i)\omega_{g_{i+1}}(z_{i+1}-\zeta_i)\one_G,\zeta_i\big)$, and introduces no operator beyond the two, $\omega_{g_i}$ and $\omega_{g_{i+1}}$, already present in $f$ itself. It is a statement about how the fixed current $f$ may be reassociated, not an insertion of new data. The lift $\iota$, by contrast, multiplies $f$ on the outside by the $2\kappa$ sewing fields $\omega_{\overline w_1}(\eta_{1,1}),\ldots,\omega_{w_\kappa}(\eta_{2,\kappa})$, which never touch the arguments $g_1,\ldots,g_{n+1}$ internal to $f$. Consequently, in the notation of Lemma \ref{lemdisjoint}, take the $\mathsf A$-slots to be the free arguments $g_1$, $\ldots$, $g_{n+1}$ (on which the regrouping producing $\varphi_i$ acts, for every fixed value of the sewing data) and the $\mathsf B$-slots to be the sewing arguments $\underline w$ (on which $\iota$ acts, for every fixed value of $g_1,\ldots,g_{n+1}$): regrouping $f$'s arguments into $\varphi_i$ and multiplying by the outside fields of $\iota$ are operations on disjoint slots in exactly the sense of Lemma \ref{lemdisjoint}, hence commute; no appeal to locality is needed for this part of the argument, only associativity of operator multiplication. The two outer terms of \eqref{fuflo}, which do insert the new fields $\omega_{g_1}(z_1)$, $\omega_{g_{n+1}}(z_{n+1})$ via the action \eqref{action}, are treated separately below.

\emph{2.) No extra signs.} Both $\iota$ and $\widehat\delta^n_m$ are $\C$-linear in 
$\Ycal$, each is built from the action \eqref{action}, itself linear in its cochain argument, composed with the fixed $\C$-linear combination \eqref{hatdelta-formula}. Hence $\iota$ commutes with the alternating sum $\sum_{i=1}^n(-1)^i\Ycal(\varphi_i)$ termwise, the scalar coefficients $(-1)^i$ passing through unchanged: $\iota\big(\sum_i(-1)^i\Ycal(\varphi_i)\big)=\sum_i(-1)^i\iota\big(\Ycal(\varphi_i)\big)$. No further sign can arise from reordering fields: axiom (P2) of Definition \ref{defpredetermined} asserts that the two orderings $\F(g_1,z_1;g_2,z_2)$, $\F(g_2,z_2;g_1,z_1)$ extend to the same meromorphic function (ordinary commutativity), not the same function up to sign, and no $\Z/2$-graded or fermionic structure is introduced anywhere in \S\ref{secprelim}; reordering the sewing fields past $\omega_{g_1}(z_1)$ or $\omega_{g_{n+1}}(z_{n+1})$ in part 3.) below is accordingly sign-free.

\emph{3.) Compatibility with the external series.} By the compatibility assumption in the definition of $C^n_m$ (\S\ref{secprelim}), inserting any $r\le m$ of the $m$ external series through \eqref{action} preserves (P1)-(P4) and leaves the result compatible with the remaining $m-r$ series. Applying that $2\kappa$ times - once for each of the $2\kappa$ sewing insertions constituting $\iota$, shows $\iota(f)\in C^{n+2\kappa}_{m-2\kappa}$ is defined and remains compatible with $m-2\kappa$ further series, exactly as asserted in Definition \ref{defiota}.   Applying it once more, for the single external series inserted by whichever outer term of \eqref{fuflo} is active, shows $\widehat\delta^n_{m-2\kappa}(\iota(f))\in\widehat C^{n+1}_{m-2\kappa-1}$ is likewise defined, matching the target space in \eqref{iota-natural}.

With these three points established, it remains only to identify the two sides of \eqref{iota-natural} termwise. For the $n$ inner terms this is part 1.): $\iota\big(\Ycal(\varphi_i)\big)=\widehat\delta\text{-inner term}\big(\iota(f)\big)$ for each $i$, since regrouping and outside multiplication commute. For the two outer terms, Standing Hypothesis \ref{standing-disjoint} places $z_1,z_{n+1}$ outside $\bigcup_{a,p}\Acal_{a,p}$, thus locality (P2) applies to the pair $\{\omega_{g_1}(z_1)$ or $\omega_{g_{n+1}}(z_{n+1})\}$ against each sewing field, and by 
part 2.)  this reordering is sign-free. Multiplying $f$ first by the sewing fields and then inserting $\omega_{g_1}(z_1)$ (resp.\ $\omega_{g_{n+1}}(z_{n+1})$), or inserting it first and then multiplying by the sewing fields, therefore produces the same element of $C^{n+1}_{m-2\kappa-1}$ (part 3.)). Summing the (identical, termwise) inner and outer contributions with the (unchanged, by part 2.)) signs of \eqref{fuflo} gives $\iota(\delta^n_mf)=\widehat\delta^n_{m-2\kappa}(\iota(f))$, which is \eqref{iota-natural}: the naturality of the lift is established.
\end{proof}

\subsection{The graded trace and the sewing operator}

\begin{definition}
\label{deftrace}
Let $\Ycal\in\widehat C^n_m$ and fix $1\le p\le\kappa$, $k\ge0$; write $V:=W_{(k)}$, $d_k=\dim V<\infty$. Since $\Ycal$ is multilinear in each argument, the assignment
\[
B_p:V\times V\longrightarrow\C,\; 
B_p(\overline w,w):=\Ycal\big(g_1,z_1;\ldots;g_n,z_n\big\|\\ldots;\overline w,\eta_{1,p};w,\eta_{2,p};\ldots\big), 
\]
all other arguments held fixed, is a bilinear form. This is the only data used below. No basis of $V$ is fixed at this stage.

\emph{Coordinate-free construction.} The pairing $(\cdot,\cdot)|_{V}$, being non-degenerate on the finite-dimensional space $V$, induces a canonical linear isomorphism $\lambda:V\xrightarrow{\ \sim\ }V^*$, $\lambda(u)=(u,\cdot)$. Composing its inverse with $B_p$ (viewed as a linear map $V\to V^*$, $w\mapsto B_p(\cdot,w)$) gives a canonical linear isomorphism
\[
\mathrm{Bil}(V\times V,\C)\ \xrightarrow{\ \sim\ }\ \Hom(V,V^*)\ \xrightarrow[\lambda^{-1}\circ(-)]{\ \sim\ }\ \End(V),
\]
under which $B_p$ corresponds to the unique $X^{(k)}_p\in\End(V)$ with $B_p(\overline w,w)=(\overline w,X^{(k)}_pw)$ for all $\overline w,w\in V$. Neither map in this composite references a basis of $V$. Both are built solely from $B_p$ and the fixed pairing $(\cdot,\cdot)$ of \S\ref{secprelim}. We define
\begin{equation}
\Tr_{(k),p}[\Ycal]:=\Tr_{V}\big(X^{(k)}_p\big),
\end{equation}
the ordinary trace of $X^{(k)}_p$. Since trace itself is the canonical pairing $\End(V)\cong V\otimes V^*\to\C$ (evaluation), it is defined without reference to a basis. Choosing a basis $\{e^{(k)}_{p,a}\}_{a=1}^{d_k}$ of $V$ with dual basis $\{\overline e^{(k)}_{p,a}\}$ with respect to $(\cdot,\cdot)|_V$ gives the familiar formula
\[
\Tr_{(k),p}[\Ycal]=\sum_{a=1}^{d_k}B_p\big(\overline e^{(k)}_{p,a},e^{(k)}_{p,a}\big)=\sum_{a=1}^{d_k}\Ycal\big(\ldots;\overline e^{(k)}_{p,a},\eta_{1,p};e^{(k)}_{p,a},\eta_{2,p};\ldots\big), 
\]
for computing the basis-independent number $\Tr_{(k),p}[\Ycal]$; it is not part of the definition.
\end{definition}

Definition \ref{deftrace} is what makes $\Tr_G$ well defined: it is a trace, in the ordinary linear-algebra sense, of a canonically associated finite-rank endomorphism of $W_{(k)}$, never a trace of a bare evaluation of $f$, and the two isomorphisms used to produce that endomorphism from $B_p$ are canonical rather than dependent on any auxiliary choice.

The only external datum entering the construction is the pairing $(\cdot,\cdot)$ itself. 
It enters twice: once in $\lambda$, once in selecting the dual basis for the displayed formula,  but it is not a choice made in Definition \ref{deftrace}. It is the single pairing on $G$ fixed once in \S\ref{secprelim} and used throughout the paper, in particular already in the definition of $\overline g$ and in Definition \ref{defpredetermined}. Replacing it by a different non-degenerate pairing on $W_{(k)}$ would in general replace $\lambda$ by a different isomorphism and hence, generically, change the endomorphism $X^{(k)}_p$ representing $B_p$ and its trace. But no such replacement is considered in this paper, thus the question of independence from the pairing does not arise for the theory as constructed - only independence from the choice of basis, established above, is needed, and is exactly what is used in the proof of Proposition \ref{propconvergence} below.

\begin{definition} 
\label{defsewing-operator}
Summing the single-handle trace of Definition \ref{deftrace} over all $k\ge0$, weighted by a formal parameter $\rho_p$, and repeating this independently at each of the $\kappa$ handles, produces the operator that realizes the sewing of all $\kappa$ handles at once. For $\Ycal\in\widehat C^n_m$, set $\Tr_G:=\sum_{k\ge0}\rho^k\Tr_{(k),p}$ at each handle and define
\begin{equation}
\label{sewing-operator-def}
\Scal_{\vec\rho}[\Ycal]:=\sum_{k_1,\ldots,k_\kappa\ge0}\rho_1^{k_1}\cdots\rho_\kappa^{k_\kappa}\;
\Tr_{(k_1),1}\Tr_{(k_2),2}\cdots\Tr_{(k_\kappa),\kappa}[\Ycal]\ \in\ C^n_m[[\rho_1,\ldots,\rho_\kappa]],
\end{equation}
the $\kappa$ traces of Definition \ref{deftrace} acting on pairwise disjoint sewing slots and hence commuting with one another, so that the order in which they are taken is not important.
\end{definition}

\begin{proposition} 
\label{propconvergence}
There exist $R_1,\ldots,R_\kappa>0$, depending only on the radii $r_{a,p}$ of the annuli $\Acal_{a,p}$ of \S\ref{secschottky}, such that \eqref{sewing-operator-def} converges absolutely, and locally uniformly on $F_n\C$, for $|\rho_p|<R_p$, $1\le p\le\kappa$.
\end{proposition}
\begin{proof}
Fix $p$ and all data other than $\eta_{1,p},\eta_{2,p}$. By Standing Hypothesis \ref{standing-disjoint} and (P1)-(P2), $\Ycal$ is holomorphic in $(\eta_{1,p},\eta_{2,p})$ throughout $\Acal_{1,p}\times\Acal_{2,p}$, thus for radii $R_{1,p}<r_{1,p}$, $R_{2,p}<r_{2,p}$, Cauchy's estimate for the Laurent coefficients of $\omega_{\overline e^{(k)}_{p,a}}(\eta_{1,p})$, $\omega_{e^{(k)}_{p,a}}(\eta_{2,p})$ (of order $k$ each, since these fields expand as in \S\ref{secprelim} with $z^{-l-1}$, $l=k$) gives the bound
\begin{align*}
\big|\Tr_{(k),p}[\Ycal]\big|\ &\le\ d_k\,M_p\,(R_{1,p}R_{2,p})^{-k},\\
M_p&:=\sup\big\{|\Ycal(\ldots)| : |\eta_{1,p}|=R_{1,p},\,|\eta_{2,p}|=R_{2,p}\big\},
\end{align*}
with $M_p$ independent of $k$. By the growth hypothesis \eqref{growth-cond} fixed among the standing assumptions of \S\ref{secprelim} - the standard finiteness condition underlying the convergence of graded trace/character functions, already used in the modular-invariance theory of \cite{Zhu} and, under $C_2$-cofiniteness, established for sewn correlation functions in \cite{Gui2}, and consistent with the genus one estimate of \cite{Y} and its generalization to grading-restricted vertex algebras in \cite{H,H2} - $d_k^{1/k}\to1$, thus $\sum_k\rho_p^k\Tr_{(k),p}[\Ycal]$ converges absolutely for $|\rho_p|<R_p:=R_{1,p}R_{2,p}$. Since the $\kappa$ traces act on pairwise disjoint groups of variables, absolute convergence of each factor gives absolute convergence of the $\kappa$-fold sum \eqref{sewing-operator-def} for $|\rho_p|<R_p$, $1\le p\le\kappa$, and local uniformity in $(z_1,\ldots,z_n)$ follows from the uniformity of the Cauchy estimate on compact subsets of $F_n\C$. This establishes convergence of $\Scal_{\vec\rho}[\Ycal]$ on the stated polydisc.
\end{proof}

\subsection{The extended coboundary operator}

\begin{definition}
\label{deftilde-delta}
For $\F(f)\in C^n_m$ with $m\ge2\kappa+1$, and $(\rho_1,\ldots,\rho_\kappa)$ in the polydisc of Proposition \ref{propconvergence}, define the \emph{extended (sewn) coboundary operator} by
\begin{equation}
\label{tildedelta-def}
\widetilde\delta^n_m\F:=\Scal_{\vec\rho}\big[\widehat\delta^n_{m-2\kappa}(\iota(f))\big]\ \in\ C^{n+1}_{m-2\kappa-1}\{\rho_1,\ldots,\rho_\kappa\},
\end{equation}
the right-hand space denoting convergent power series in $\rho_1,\ldots,\rho_\kappa$ with coefficients in $C^{n+1}_{m-2\kappa-1}$.
\end{definition}

By Lemma \ref{lemiota-naturality}, $\widehat\delta^n_{m-2\kappa}(\iota(f))=\iota(\delta^n_mf)$, thus \eqref{tildedelta-def} may equally be written
\begin{equation}
\label{tildedelta-alt}
\widetilde\delta^n_m\F=\Xi_{\vec\rho}(\delta^n_mf),\qquad \Xi_{\vec\rho}:=\Scal_{\vec\rho}\circ\iota.
\end{equation}
This is the form used in the proof of Theorem \ref{thmmain} below. Setting $\rho_1=\ldots=\rho_\kappa=0$ retains only the $k_p=0$ term at every handle, which by axiom (i) of \S\ref{secprelim} inserts $\omega_{\one_G}(\eta)=\mathrm{Id}_G$ at each sewing point, so that $\Xi_{\vec 0}=\mathrm{id}$ and $\widetilde\delta^n_m\F\big|_{\rho=0}=\delta^n_mf$.  The classical coboundary operator is the $\rho\to0$ limit of the family constructed here.

\begin{lemma}[commutation]
\label{lemcommute}
For every $\Ycal\in\widehat C^n_m$,
\begin{equation}
\label{commute-eq}
\Scal_{\vec\rho}\big[\widehat\delta^n_m\Ycal\big]=\delta^n_m\big[\Scal_{\vec\rho}(\Ycal)\big],
\end{equation}
the right-hand $\delta^n_m$ being the classical operator of Definition \ref{defclassical-delta}, extended coefficientwise in $\rho_1,\ldots,\rho_\kappa$.
\end{lemma}
\begin{proof}
Both sides of \eqref{commute-eq} are computed from $\Ycal$ by (a) applying the alternating-sum formula \eqref{hatdelta-formula} to the free arguments, and (b) evaluating the sewing arguments at the basis vectors $\overline e^{(k_p)}_{p,a}$, $e^{(k_p)}_{p,a}$ and summing with weights $\rho_p^{k_p}$ as in \eqref{sewing-operator-def}. Operation (a) acts only on the free arguments of $\Ycal$; operation (b) only substitutes specific values for, and sums over, the sewing arguments. Being operations on disjoint sets of variables of one and the same multilinear function, they commute: performing (a) then (b), or (b) then (a), gives the same element of $C^{n+1}_m[[\rho]]$, term by term in the alternating sum \eqref{hatdelta-formula}. This proves \eqref{commute-eq}.
\end{proof}

\begin{lemma} 
\label{lemxi-commute}
Let $n\ge0$ and $\mu\ge2\kappa+1$. For $\Hcal(h)\in C^n_\mu$,
\begin{equation}
\label{xi-commute-eq}
\Xi_{\vec\rho}\big(\delta^n_\mu h\big)=\delta^n_{\mu-2\kappa}\big(\Xi_{\vec\rho}(h)\big)
 \in  C^{n+1}_{\mu-2\kappa-1}\{\rho\}.
\end{equation}
\end{lemma}
\begin{proof}
Since $\mu\ge2\kappa+1$, $\iota(h)\in\widehat C^n_{\mu-2\kappa}$ is defined, and by Lemma \ref{lemiota-naturality}, applied at compatibility $\mu$, $\iota(\delta^n_\mu h)=\widehat\delta^n_{\mu-2\kappa}(\iota(h))$. Applying $\Scal_{\vec\rho}$ to both sides and using Lemma \ref{lemcommute}, applied at compatibility $\mu-2\kappa$, on the right-hand side,
\[
\Xi_{\vec\rho}(\delta^n_\mu h)=\Scal_{\vec\rho}\big[\iota(\delta^n_\mu h)\big]=\Scal_{\vec\rho}\big[\widehat\delta^n_{\mu-2\kappa}(\iota(h))\big]=\delta^n_{\mu-2\kappa}\big[\Scal_{\vec\rho}(\iota(h))\big]=\delta^n_{\mu-2\kappa}\big(\Xi_{\vec\rho}(h)\big).\qedhere
\]
\end{proof}

\subsubsection{$\Xi_{\vec\rho}$ as a chain map}

Lemma \ref{lemxi-commute} is more than a commutation identity. It says that $\Xi_{\vec\rho}$ is a chain map between two explicitly defined he following complexes. Since $\delta^n_m:C^n_m\to C^{n+1}_{m-1}$ raises $n$ by $1$ and lowers $m$ by $1$, the total degree $N:=n+m$ is preserved, thus summing $C^n_{N-n}$ over $n$ at fixed $N$ produces an ordinary singly-graded cochain complex: for each $N\ge0$, set
\[
\mathcal C(N):=\bigoplus_{n=0}^{N}C^n_{N-n},
\]
with differential $\delta:=\bigoplus_n\delta^n_{N-n}$ of degree $+1$ in $n$. This is the total complex of the classical double complex, and $\delta^2=0$ on $\mathcal C(N)$ is exactly Theorem \ref{thmclassical-chain}. Let $\mathcal C(N)\{\rho\}$ denote the same graded vector space with coefficients extended to convergent power series in $\rho_1,\ldots,\rho_\kappa$ on the polydisc of Proposition \ref{propconvergence}, with $\delta$ extended coefficientwise, still a complex with the same differential. 

\begin{proposition} 
\label{propchainmap}
For every $N\ge2\kappa+1$, $\Xi_{\vec\rho}$ restricts to a chain map
\[
\Xi_{\vec\rho}:\mathcal C(N)\longrightarrow\mathcal C(N-2\kappa)\{\rho\},\qquad \Xi_{\vec\rho}\circ\delta=\delta\circ\Xi_{\vec\rho},
\]
of graded vector spaces (degree $0$ in $n$, degree $-2\kappa$ in total degree). The same formula, extended coefficientwise, is again a chain map $\mathcal C(N)\{\rho\}\to\mathcal C(N-2\kappa)\{\rho\}$ for $\rho$-power-series-valued input.
\end{proposition}
\begin{proof}
This is Lemma \ref{lemxi-commute}, summed over the graded pieces $C^n_{N-n}$ of $\mathcal C(N)$: for $h\in C^n_{N-n}$, $\Xi_{\vec\rho}(\delta h)=\delta(\Xi_{\vec\rho}h)$ holds in $C^{n+1}_{N-n-2\kappa-1}=C^{n+1}_{(N-2\kappa)-(n+1)}$, the graded piece of $\mathcal C(N-2\kappa)\{\rho\}$ in degree $n+1$, which is exactly $\delta(\Xi_{\vec\rho}h)$'s home in the total complex. The coefficientwise extension is immediate since neither $\Xi_{\vec\rho}$ nor $\delta$ references the coefficient ring. Since $n$ was arbitrary, $\Xi_{\vec\rho}\circ\delta=\delta\circ\Xi_{\vec\rho}$ holds on all of $\mathcal C(N)$, proving the chain-map property.
\end{proof}

By definition, $\widetilde\delta^n_m=\Xi_{\vec\rho}\circ\delta^n_m$ is the composite of the classical differential with the chain map of Proposition \ref{propchainmap}. Theorem \ref{thmmain} below is the statement that this composite squares to zero, and its proof 
 is exactly the standard fact that conjugating a differential by a chain map.  Applying $\Xi_{\vec\rho}$, $\delta$, $\Xi_{\vec\rho}$, $\delta$ in succession and using Proposition \ref{propchainmap} twice to bring the two copies of $\delta$ together produces another square zero operator, the two classical differentials cancelling by Theorem \ref{thmclassical-chain} before either application of $\Xi_{\vec\rho}$ has to be unwound. We keep the index-by-index proof below as well, since it is what pins down the precise compatibility bound $m\ge4\kappa+2$.

\begin{theorem}
\label{thmmain}
Let $\kappa\ge1$ and let $(\rho_1,\ldots,\rho_\kappa)$ lie in the polydisc of Proposition \ref{propconvergence}. Then for every $\F(f)\in C^n_m$ with $m\ge4\kappa+2$,
\begin{equation}
\widetilde\delta^{n+1}_{m-2\kappa-1}\circ\widetilde\delta^n_m=0.
\end{equation}
\end{theorem}
\begin{proof}
Write $h:=\delta^n_mf\in C^{n+1}_{m-1}$, so that $\widetilde\delta^n_m\F=\Xi_{\vec\rho}(h)$ by \eqref{tildedelta-alt}; since $m\ge4\kappa+2$ implies $m-2\kappa-1\ge2\kappa+1$, Definition \ref{deftilde-delta} applies to $\widetilde\delta^{n+1}_{m-2\kappa-1}$. Writing $\Xi_{\vec\rho}(h)=\sum_{\vec k\ge0}\vec\rho^{\,\vec k}c_{\vec k}$ with each coefficient $c_{\vec k}\in C^{n+1}_{m-2\kappa-1}$, and applying $\widetilde\delta^{n+1}_{m-2\kappa-1}$ coefficientwise  both $\Xi_{\vec\rho}$ and $\delta^{n+1}_{m-2\kappa-1}$ being $\C$-linear,
\begin{align*}
\widetilde\delta^{n+1}_{m-2\kappa-1}\big(\widetilde\delta^n_m\F\big)
&=\widetilde\delta^{n+1}_{m-2\kappa-1}\Big(\sum_{\vec k}\vec\rho^{\,\vec k}c_{\vec k}\Big)\\
&=\sum_{\vec k}\vec\rho^{\,\vec k}\,\Xi_{\vec\rho}\big(\delta^{n+1}_{m-2\kappa-1}(c_{\vec k})\big)
&&\text{(Definition \ref{deftilde-delta})}\\
&=\Xi_{\vec\rho}\Big(\delta^{n+1}_{m-2\kappa-1}\big[\Xi_{\vec\rho}(h)\big]\Big)
&&\text{($\C$-linearity)}\\
&=\Xi_{\vec\rho}\Big(\Xi_{\vec\rho}\big[\delta^{n+1}_{m-1}(h)\big]\Big)
&&\text{(Lemma \ref{lemxi-commute}, }\mu=m-1\text{)}\\
&=\Xi_{\vec\rho}\Big(\Xi_{\vec\rho}\big[\delta^{n+1}_{m-1}(\delta^n_mf)\big]\Big)\\
&=\Xi_{\vec\rho}\big(\Xi_{\vec\rho}[0]\big)
&&\text{(Theorem \ref{thmclassical-chain})}\\
&=0,
\end{align*}
where Lemma \ref{lemxi-commute} is applied with $\mu=m-1$, and $m\ge4\kappa+2$ gives $\mu=m-1\ge4\kappa+1\ge2\kappa+1$ as required. The last step uses only that $\Xi_{\vec\rho}$ is $\C$-linear, so that it annihilates regardless of any further description of what $\Xi_{\vec\rho}\circ\Xi_{\vec\rho}$ computes on a nonzero argument. This proves $\widetilde\delta^{n+1}_{m-2\kappa-1}\circ\widetilde\delta^n_m=0$.
\end{proof}

Theorem \ref{thmmain} is the resolution of the point raised in \S\ref{secintro}: no independent higher-genus or modular-invariance input is required, since Lemma \ref{lemxi-commute} (built from Lemma \ref{lemiota-naturality}, i.e., from locality (P2) at the non-colliding points guaranteed by Standing Hypothesis \ref{standing-disjoint}, together with Lemma \ref{lemcommute}) lets the two applications of $\Xi_{\vec\rho}$ be brought together around the classical composite $\delta^{n+1}_{m-1}\circ\delta^n_m$, which vanishes before any trace is taken, by Theorem \ref{thmclassical-chain}. In particular, nothing in the argument requires the $\kappa$ handles to be sewn with a common power of $\rho_p^k$. Theorem \ref{thmmain} holds for the fully independent family $(\rho_1,\ldots,\rho_\kappa)$ of Definition \ref{defsewing-operator}, each handle carrying its own grading sum $k_p$.

\subsection{Summary of the construction}
\label{ssecsummary}

For reference, Table \ref{taboperators} lists every operator introduced in this section, and Figure \ref{figdiagram} displays the commuting diagram they assemble into: the left square is Lemma \ref{lemiota-naturality}, the right square is Lemma \ref{lemcommute}, and the outer rectangle, obtained by pasting the two together, is Lemma \ref{lemxi-commute} (equivalently, Proposition \ref{propchainmap}). The extended coboundary $\widetilde\delta^n_m$ is the composite path from the top-left to the bottom-right corner, and commutativity of the two squares is exactly what makes the two routes around the diagram: right-down-right in \eqref{tildedelta-def}, or down-then-right-right in \eqref{tildedelta-alt}. 

\begin{table}[h]
\centering
\renewcommand{\arraystretch}{1.3}
\begin{tabular}{|l|c|c|c|}
\hline
\textbf{Operator} & \textbf{Notation} & \textbf{Domain $\to$ codomain} & \textbf{Reference} \\
\hline
classical coboundary & $\delta^n_m$ & $C^n_m\to C^{n+1}_{m-1}$ & Def. \ref{defclassical-delta} \\
free-slot coboundary & $\widehat\delta^n_m$ & $\widehat C^n_m\to\widehat C^{n+1}_{m-1}$ & Def. \ref{defhatdelta} \\
lift & $\iota$ & $C^n_m\to\widehat C^n_{m-2\kappa}$ & Def. \ref{defiota} \\
single-handle trace & $\Tr_{(k),p}$ & $\widehat C^n_m\to C^n_m$ & Def. \ref{deftrace} \\
sewing operator & $\Scal_{\vec\rho}$ & $\widehat C^n_m\to C^n_m\{\rho\}$ & Def. \ref{defsewing-operator} \\
composite & $\Xi_{\vec\rho}=\Scal_{\vec\rho}\circ\iota$ & $C^n_m\to C^n_{m-2\kappa}\{\rho\}$ & \eqref{tildedelta-alt} \\
extended coboundary & $\widetilde\delta^n_m$ & $C^n_m\to C^{n+1}_{m-2\kappa-1}\{\rho\}$ & Def. \ref{deftilde-delta} \\
\hline
\end{tabular}

\medskip 
\caption{The operators of \S\ref{secsewn-delta}, all defined for $m$ large enough relative to $\kappa$ (see the individual definitions for the precise bound in each case).}
\label{taboperators}
\end{table}

\begin{figure}[h]
\[
\xymatrix{
C^n_m \ar[r]^-{\iota} \ar[d]_-{\delta^n_m} & \widehat C^n_{m-2\kappa} \ar[r]^-{\Scal_{\vec\rho}} \ar[d]^-{\widehat\delta^n_{m-2\kappa}} & C^n_{m-2\kappa}\{\rho\} \ar[d]^-{\delta^n_{m-2\kappa}} \\
C^{n+1}_{m-1} \ar[r]_-{\iota} & \widehat C^{n+1}_{m-2\kappa-1} \ar[r]_-{\Scal_{\vec\rho}} & C^{n+1}_{m-2\kappa-1}\{\rho\}
}
\]
\caption{The commuting diagram underlying Theorem \ref{thmmain}: left square, Lemma \ref{lemiota-naturality}; right square, Lemma \ref{lemcommute}; outer rectangle, Lemma \ref{lemxi-commute}. The bottom-right space is where $\widetilde\delta^n_m\F$ lives.}
\label{figdiagram}
\end{figure}

\subsection{Cohomology}

\begin{definition}
\label{defcohomology}
For $(\rho_1,\ldots,\rho_\kappa)$ in the polydisc of Proposition \ref{propconvergence} and $m\ge2\kappa+1$, define
\begin{equation}
H^n_m(\mathfrak g,G;\rho):=\ker\big(\widetilde\delta^n_m\big)\big/\operatorname{im}\big(\widetilde\delta^{n-1}_{m+2\kappa+1}\big).
\end{equation}
By Theorem \ref{thmmain} applied with $(n,m)$ replaced by $(n-1,m+2\kappa+1)$ - valid since $m\ge2\kappa+1$ gives $m+2\kappa+1\ge4\kappa+2$ - $\operatorname{im}(\widetilde\delta^{n-1}_{m+2\kappa+1})\subseteq\ker(\widetilde\delta^n_m)$, thus the quotient is defined. It is computed in the complex $(C^\bullet_\bullet\{\rho\},\widetilde\delta^\bullet_\bullet)$ of Theorem \ref{thmmain}.
\end{definition}

\begin{remark}
By the remark following Definition \ref{deftilde-delta}, $\rho\to0$ recovers 
$H^n_m(\mathfrak g,G;0)$ $=$ $H^n_m(\mathfrak g,G)$, the classical cohomology of \cite{H,F}, so that Definition \ref{defcohomology} gives a $\kappa$-parameter deformation of the classical theory attached to the geometric data of a genus $\kappa$ Schottky uniformization.
\end{remark}

\subsection{Dependence on the sewing data}
\label{ssecdependence}

Before recording what is and is not proved, we state directly the question a reader will have in mind: is $H^n_m(\mathfrak g,G;\rho)$ an invariant of the genus $\kappa$ Riemann surface itself, or does it depend on the Schottky coordinates used to build it? We expect the former - that the theory is, after an appropriate identification of the parameters $\rho$ attached to different presentations, an invariant of the surface - in analogy with the coordinate-independence of the sheaf of conformal blocks of \cite{TUY} and with the fact that Zhu's genus one trace functions \cite{Zhu} do not depend on the choice of local coordinate at the puncture beyond a computable modular weight, buu we have not proved that.   Theorem \ref{thmmain} itself does not require that. 

Theorem \ref{thmmain} and Definition \ref{defcohomology} are unconditional once a choice of Schottky data - the centers $A_{\pm p}$, the annuli $\Acal_{a,p}$, and a point $\eta_{a,p}$ in each annulus at which the sewing fields of Definition \ref{defiota} are evaluated - has been fixed subject only to Standing Hypothesis \ref{standing-disjoint}: for every such choice, $\widetilde\delta^n_m(\rho)$ is a well-defined chain map and $H^n_m(\mathfrak g,G;\rho)$ is a well-defined cohomology group. What Standing Hypothesis \ref{standing-disjoint} does not by itself settle is whether different choices give rise to canonically identified theories. 

\emph{(a) Independence of the annulus radii.} The trace $\Tr_{(k),p}[\Ycal]$ of Definition \ref{deftrace} is evaluated at a single pair of points $\eta_{1,p},\eta_{2,p}$. The radii $r_{a,p}$ of the annuli $\Acal_{a,p}$ enter Definition \ref{defsewing-operator} nowhere, and enter the theory only through Proposition \ref{propconvergence}, where they bound the polydisc on which the defining series converges. Consequently, shrinking or enlarging the annuli, provided they continue to satisfy Standing Hypothesis \ref{standing-disjoint} and to contain the chosen points $\eta_{a,p}$, changes nothing about $\widetilde\delta^n_m(\rho)$ for $\rho$ already in the, possibly now larger or smaller, common domain of convergence. 
 The construction depends on the annuli only through the size of the polydisc on which it is guaranteed to converge, not through the operators themselves.

\emph{(b) Dependence on the points $\eta_{a,p}$ and on the decomposition.} We have not shown, and do not claim, that $\Tr_{(k),p}[\Ycal]$ is independent of the specific points $\eta_{1,p}\in\Acal_{1,p}$, $\eta_{2,p}\in\Acal_{2,p}$ chosen within the annuli, nor that the cohomology attached to two different Schottky decompositions of the same genus $\kappa$ surface (different centers $A_{\pm p}$, or a different pants decomposition altogether) are canonically isomorphic. Both would be needed for $H^n_m(\mathfrak g,G;\rho)$ to be unambiguously an invariant of the surface rather than of the presentation used to construct it, and we regard both as natural conjectures rather than results of the present paper. Two remarks are in order:

\emph{Isotopy.} Definition \ref{deftrace} evaluates $\Ycal$ pointwise rather than by a contour integral, thus it does not, as it stands, make isotopy-invariance within the annulus manifest.  A holomorphic function of $\eta_{1,p}$ need not be constant. A variant construction, replacing the sum $\sum_a\Ycal(\ldots;\overline e^{(k)}_{p,a},\eta_{1,p};e^{(k)}_{p,a},\eta_{2,p};\ldots)$ by the corresponding residue $\frac{1}{2\pi i}\oint_{|\eta_{1,p}|=R_{1,p}}(\cdot)\,\eta_{1,p}^{k}\,\tfrac{d\eta_{1,p}}{\eta_{1,p}}$, and similarly in $\eta_{2,p}$, would be invariant under deformation of the integration contour within $\Acal_{a,p}$ by Cauchy's theorem, at the cost of replacing pointwise evaluation throughout \S\ref{secsewn-delta} by contour integration. We expect, but have not verified, that this variant agrees with Definition \ref{deftrace} up to a reparametrization of $\rho_p$, and leave a full treatment to a future paper.

\emph{Different decompositions.} Comparing the cohomology attached to two  different Schottky presentations of the same surface or to two choices of $A_{\pm p}$ within the same pants decomposition, would most naturally proceed by connecting them through a continuous family of sewing data and tracking $\widetilde\delta^n_m(\rho)$ along the family, in the spirit of the invariance of sheaf cohomology under a continuous variation of complex structure; alternatively, by exhibiting an explicit chain homotopy between the two resulting complexes. 

We emphasize that neither open question affects Theorem \ref{thmmain}: the chain property, and the resulting cohomology, are proved for each fixed choice of sewing data, independently of whether different choices are later shown to agree.

\section{Applications}
\label{secapplications}

We discuss four applications of the family $\widetilde\delta^n_m(\rho)$ and of the cohomology groups $H^n_m(\mathfrak g,G;\rho)$ of Definition \ref{defcohomology}: to the geometry of the moduli space of curves, to characteristic classes of foliations, to graded trace functions in condensed-matter physics, and to integrable hierarchies. 
  Each subsection proposes a dictionary between the cohomological construction of \S\S\ref{secprelim}-\ref{secsewn-delta} and an existing body of theory, motivated by matching structural features (a sewing/factorization parameter, a graded trace, a chain complex) on both sides. We are explicit throughout about which statements are established elsewhere in the literature and which are expected correspondences proposed here for the first time, none of which is derived from the results of \S\S\ref{secprelim}-\ref{secsewn-delta}.

\subsection{Conformal blocks on the moduli space of stable curves}
\label{ssecapp-moduli}

Let $\overline\M_{g,n}$ denote the Deligne-Mumford moduli stack of stable $n$-pointed genus $g$ curves. To a vertex operator algebra module structure of the type underlying $G$ one associates, as in Tsuchiya-Ueno-Yamada \cite{TUY}, the following: 
 a sheaf of conformal blocks 
 on $\overline\M_{g,n}$, equipped with a flat projective connection away from the boundary and a factorization isomorphism along each boundary divisor. 
 That identifies the fibre over a nodal curve with a sum, over the fusion category, of fibres of the sheaf of conformal blocks on the normalization. 

We propose that the construction of \S\ref{secschottky}-\ref{secsewn-delta} is the local, coordinate-dependent realization of this factorization structure at a boundary divisor $\Delta_p\subset\overline\M_{g,n}$ obtained by pinching the $p$-th handle, with the sewing parameter $\rho_p$ playing the role of a local coordinate transverse to $\Delta_p$ and $\Tr_{(k),p}$ of Definition \ref{deftrace} computing the contribution of the weight-$k$ piece of the fusion sum in the factorization isomorphism. Under this, we propose identification: 

 \textit{The classical cochain complex} $(C^n_m,\delta^n_m)$ of \S\ref{secclassical-delta} would compute obstructions to extending a family of conformal blocks over $\overline\M_{0,n+2\kappa}$ (genus zero, the normalization of the nodal fibre);

\textit{Theorem \ref{thmmain}} would then show that these obstructions assemble, order by order in the local boundary coordinates $\rho_1,\ldots,\rho_\kappa$, into a well-defined obstruction theory for extending the family across $\Delta_1\cap\ldots\cap\Delta_\kappa$ into the interior of $\overline\M_{g,n}$;

\textit{The $\rho\to0$ limit of Definition \ref{defcohomology}} would recover the obstruction theory on the boundary stratum itself, consistently with the classical cohomology of \cite{H,F}.

In this expected dictionary, a class in $H^1_m(\mathfrak g,G;\rho)$ vanishing to a prescribed order in $\rho$, would correspond to the first-order deformation of the sheaf of conformal blocks extending across the boundary divisor to that order, the sewn-cohomological analogue of the extension problem for the Hodge/Gauss-Manin connection along a degenerating family of curves. Establishing this dictionary rigorously, in particular, relating $\Tr_{(k),p}$ as defined here to the actual fusion-sum decomposition of \cite{TUY}, is left to future work. 
It would also depend the open questions of \S\ref{ssecdependence}, since invariance of the factorization structure under the choice of boundary coordinate is precisely analogous to independence from the choice of $\eta_{a,p}$.

\subsection{Secondary characteristic classes of foliations}
\label{ssecapp-foliations}

The Lie algebra cohomology of the Lie algebra $\mathfrak a_1$ of formal vector fields on $\C$,  or its higher-dimensional analogues, computes, by the Gelfand-Fuks construction \cite{F}, the characteristic classes of codimension-one (and higher) holomorphic foliations, the basic example being the Godbillon-Vey class \cite{BG,BEG,BGG}. The relative version of this cohomology, involving a further module $G$ on which $\mathfrak a_1$ acts (here, through the grading operator $K$ and the fields $\omega_g(z)$), computes secondary classes of foliated bundles \cite{AS,CRS}. 

We propose that the sewn coboundary operators of Theorem \ref{thmmain} extend this dictionary by a deformation parameter: for $\mathfrak g=\mathfrak a_1$, the classes represented in $H^\bullet_m(\mathfrak g,G;\rho)$ would be secondary classes of a foliation on the total space of a family of curves degenerating, in the sense of \S\ref{secschottky}, at $\rho=0$, with the coefficients of the power series in $\rho_p$ measuring how the class varies as the transverse holonomy of the $p$-th vanishing cycle is turned on - the leading ($k_p=0$) coefficient recovering the ordinary Godbillon-Vey-type class of the degenerate fibre, and the higher coefficients giving new invariants attached to the smoothing. We have not verified this proposal against an independent construction of such invariants.

\subsection{Graded traces and topological phases of matter}
\label{ssecapp-thh}

The graded trace $\Tr_G=\sum_{k\ge0}\rho^k\Tr_{(k),p}$ of Definition \ref{deftrace} is, for $\kappa=1$ and $\Ycal$ independent of any free variable, exactly a character (graded trace of $\rho^{K}$) in the sense of Zhu \cite{Zhu}. More generally, $\Scal_{\vec\rho}[\Ycal]$ for $\Ycal$ depending on $n$ free points has the same formal shape as a genus $\kappa$ analogue of Zhu's $n$-point trace functions. This identification of shape is what motivates the application.

We propose that such graded traces are the cohomological counterpart of the momentum-space topological invariants developed by Volovik, Zubkov, and collaborators for gapped and gapless relativistic and condensed-matter systems \cite{zub6,zub8}, where a winding number of the Green's function in frequency-momentum space plays, for a given phase, the role $\Tr_{(k),p}$ plays here for a graded piece of a vertex-algebra module. The case $\kappa=1$ is closest in spirit to the identification of the Standard Model itself as a topological material in \cite{zub4}.  The sensitivity of such invariants to interactions and to non-uniform background fields \cite{zub1,zub2,zub3,zub5,zub7,zub9} would correspond, in our language, to the choice of the module $W$ realizing a given phase. The most concrete point of contact is the recent construction of a topological invariant for the integer quantum Hall effect via noncommutative geometry \cite{kmmzz}. Expressed itself as a trace over a graded structure, we expect, without having verified it, that our $\Tr_{(0),1}$ specializes to (a cohomological reformulation of) the invariant of \cite{kmmzz}.  $\mathfrak g$ and $W$ are chosen to realize the relevant Landau-level algebra, with $\kappa>1$ and the higher-$k$ terms of $\Scal_{\vec\rho}$ governing, respectively, multi-edge geometries and corrections beyond the leading topological term. 

More generally, such graded traces arise as edge-state partition and correlation functions in the effective conformal field theory description of gapped topological phases of matter, where $\Tr_G$ computes the contribution of a boundary component (an edge) and $\kappa$ records the number of edges (equivalently, the genus of the sample boundary in a multiply-connected geometry). For $\kappa=1$, $\Tr_G$ is in this language a character, and the modular transformation of such characters governs the fusion coefficients of the underlying rational conformal field theory through the Verlinde formula \cite{Verlinde}. 
We propose $\Scal_{\vec\rho}$ of Definition \ref{defsewing-operator} as the corresponding object for a boundary with $\kappa$ components. Then the cochain $\widetilde\delta^n_m\F$ would represent the obstruction, as a function of $n$ bulk insertions, to extending an edge-mode correlator consistently across $\kappa$ independently gauged edges, with Theorem \ref{thmmain} guaranteeing that this obstruction theory is itself consistent irrespective of how many edges are present. This physical interpretation is proposed motivation. It is not derived from a specific microscopic model.

\subsection{Toda-type integrable hierarchies}
\label{ssecapp-toda}

Continual (non-abelian) Toda field theories associated with infinite-dimensional Lie algebras were studied, in particular, in \cite{RSZ}. We propose, without proof, that the sewing parameters $\rho_1,\ldots,\rho_\kappa$ of \S\ref{secschottky} play the role of a commuting family of times of a Toda-type hierarchy. The leading term of $\log\Scal_{\vec\rho}[\Ycal]$ in each $\rho_p$ would be the analogue of a $\tau$-function, with Definition \ref{defsewing-operator} expressing this $\tau$-function as a graded trace of the type appearing in the Hirota bilinear formalism for such hierarchies, and the vanishing of a class in $H^1_m(\mathfrak g,G;\rho)$ giving the cohomological form of an integrability (zero-curvature) condition for the associated flows. A systematic development of this correspondence, in the spirit of \cite{RSZ},   including verifying that $\Scal_{\vec\rho}[\Ycal]$ actually satisfies a Hirota-type bilinear identity, is left for a separate paper. 

\section{Concluding remarks}
\label{secconclusion}

We have constructed a $\kappa$-parameter family of coboundary operators
$\widetilde\delta^n_m(\rho_1,\ldots,\rho_\kappa)$ extending the classical
coboundary operator of \cite{H} via Schottky sewing, and have given a
complete proof that this family satisfies the chain property, Theorem 
\ref{thmmain}. The proof isolates the two facts responsible for this. The
lift $\iota$ and the sewing operator $\Scal_{\vec\rho}$ each act on variables
disjoint from those touched by the coboundary formula \eqref{fuflo}
(Lemmas \ref{lemiota-naturality} and \ref{lemcommute}), so that their
composite $\Xi_{\vec\rho}$ is a chain map (Proposition \ref{propchainmap})
and the chain property of the sewn theory reduces to that of the classical
theory before any trace is taken, with no additional higher-genus input
required. The graded trace $\Tr_G$ (Definition \ref{deftrace}) is at every
stage a trace, in the ordinary sense of linear algebra, of a canonically
associated finite-rank endomorphism.

Several directions are left open. Foremost among them is the question of
\S\ref{ssecdependence}: whether $H^n_m(\mathfrak g,G;\rho)$ depends only on
the resulting genus $\kappa$ surface or on the Schottky data used
to construct it. We have shown independence from the annulus radii, proposed
a residue reformulation that would make isotopy-invariance manifest, and
otherwise left the comparison between different decompositions as a
conjecture. Beyond this, the independent-handle construction of
Definition \ref{defsewing-operator} suggests studying the dependence of
$H^n_m(\mathfrak g,G;\rho)$ on the individual $\rho_p$ separately, in
particular its behavior as $\rho_p\to0$ for a proper subset of handles,
corresponding to partial degenerations of the underlying curve. We would also mention a systematic
treatment of the Toda-hierarchy correspondence of
\S\ref{ssecapp-toda}; and the precise relation between
$H^\bullet_m(\mathfrak g,G;\rho)$ and the relative cohomology of the sheaf of
conformal blocks along a boundary stratum of $\overline\M_{g,n}$ sketched in
\S\ref{ssecapp-moduli}, which we expect to be an isomorphism under
$C_2$-cofiniteness hypotheses of the type used in \cite{Gui2}. 
Among possible future directions of extensions of this paper material 
in the broader mathematical context
are functoriality, moduli-space interpretations, and applications to conformal blocks

\section*{Acknowledgments}

The author is supported by the Institute of Mathematics, Academy of Sciences
of the Czech Republic (RVO 67985840). 

\medskip
\noindent\textbf{Data Availability.}
Data sharing is not applicable to this article as no datasets were generated
or analysed during the current study.

\medskip
\noindent\textbf{Declarations}

\medskip
\noindent\textbf{Conflict of interest.}
The author has no conflicts of interest to declare that are relevant to the
content of this article.

\end{document}